\PassOptionsToPackage{off}{crop}

\documentclass[preprint,ejsv2,,noshowframe]{imsart}

\usepackage{fancyhdr}
\makeatletter
\AtBeginDocument{%
  \@ifundefined{CROP@offtrue}{}{\CROP@offtrue}%
  \@ifundefined{CROP@infofalse}{}{\CROP@infofalse}%
}
\makeatother

\usepackage[authoryear,round]{natbib}
\usepackage{graphicx} 
\usepackage{amsmath,amsthm,bm,xcolor}

\newcommand{\e}{\mathrm{e}}
\providecommand{\arrowvert}{\vert}

\providecommand{\skakko}[1]{\left(#1\right)}

\numberwithin{equation}{section}
\theoremstyle{plain}
\newtheorem{thm}{Theorem}[section]
\newtheorem{lem}[thm]{Lemma}
\newtheorem{cor}[thm]{Corollary}

\theoremstyle{definition}
\newtheorem{rem}{Remark} 

\startlocaldefs

\theoremstyle{definition}

\theoremstyle{remark}

\endlocaldefs

\begin{document}
\begin{frontmatter}
\title{On spectral clustering under non-isotropic Gaussian mixture models}
\runtitle{Spectral clustering}

\begin{aug}
\author[A]{\fnms{Kohei}~\snm{Kawamoto}%
  \ead[label=e1]{kawamoto.kohei.532@s.kyushu-u.ac.jp}%
  \orcid{0009-0001-9424-0018}}
\author[B]{\fnms{Yuichi}~\snm{Goto}%
  \ead[label=e2]{yuichi.goto@math.kyushu-u.ac.jp}%
  \orcid{0000-0002-7556-2572}}
\and
\author[B]{\fnms{Koji}~\snm{Tsukuda}%
  \ead[label=e3]{tsukuda@math.kyushu-u.ac.jp}%
  \orcid{0000-0002-3973-0460}}
\address[A]{Joint Graduate School of Mathematics for Innovation, Kyushu University,
744 Motooka, Fukuoka, Fukuoka 819-0395, Japan\printead[presep={,\ }]{e1}}

\address[B]{Faculty of Mathematics, Kyushu University,
744 Motooka, Fukuoka, Fukuoka 819-0395, Japan\printead[presep={,\ }]{e2,e3}}

\runauthor{K. Kawamoto et al.}
\end{aug}
\begin{abstract}
We evaluate the misclustering probability of a spectral clustering algorithm under a Gaussian mixture model with a general covariance structure. 
The algorithm partitions the data into two groups based on the sign of the  first principal component score.
As a corollary of the main result, the clustering procedure is shown to be consistent in a high-dimensional regime.
\end{abstract}

\begin{keyword}[class=MSC2020]
\kwd{62H30}
\kwd{62H25}
\kwd{62H12}
\end{keyword}

\begin{keyword}
\kwd{clustering}
\kwd{high-dimensional regime}
\kwd{misclustering probability}
\end{keyword}

\end{frontmatter}

\section{Introduction}
Spectral clustering algorithm is a widely used method that applies a clustering algorithm, such as $k$-means, after dimension reduction via spectral decomposition of a certain matrix, including an adjacency matrix, a graph Laplacian derived from a similarity matrix, or a covariance matrix.
Compared with other clustering approaches, spectral clustering algorithm has the advantage of a lower computational burden.
Unlike iterative methods, spectral clustering algorithm does not require tuning parameters.

Considerable attention has been devoted to spectral clustering algorithm based on 
$k$-means for more than two clusters; see, for example,
\cite{lzz21}, \citet[Section~3.1]{afw22}, \cite{zz24}, and references therein.
In contrast, for binary clustering problems, sign-based spectral clustering algorithm methods have also been extensively investigated.
\cite{cz18}, \citet[Section~3.2]{afw22}, and \citet[Section~4.7.1]{v18} considered clustering based on the sign of the first right singular vector, the leading eigenvector of the hollowed Gram matrix, and the first principal component score, respectively. These studies assume the Gaussian mixture model with isotropic innovations.
Under a heteroskedastic two-component Gaussian mixture model with independent coordinates and coordinate-wise variances, \citet{chz22} studied clustering based on the sign of the first principal component score, leveraging sharp
non-asymptotic concentration bounds for Wishart-type matrices.
More recently, \cite{kgt25} analyzed this sign-based spectral clustering procedure under an allometric extension relationship, in which the first principal directions of the two clusters are aligned with the direction of the difference between their mean vectors.

As stated, existing theoretical studies assume certain restrictions on the covariance structure such as isotropy or allometric extension model.
The leading principal component of their mixture models other than the model of \cite{chz22} aligns with the mean difference direction.
However, these model assumptions often do not fit real-world data.
In this paper, we derive a non-asymptotic bound on the misclustering probability of the spectral clustering algorithm method for a two-component Gaussian mixture model without such restrictions and, as a by-product, sufficient conditions ensuring consistency of spectral clustering in certain high-dimensional regimes.


\section{Model and spectral clustering algorithm}
Let $\theta$ be a Rademacher random variable that takes a value in $\{-1,1\}$ with equal probability. Denote 
$\bm{g}^{(-1)}\sim \mathcal{N}_n(\bm{0}_n,\bm{\Sigma}_{-1})$ and $\bm{g}^{(1)}\sim \mathcal{N}_n(\bm{0}_n,\bm{\Sigma}_1)$ for positive-definite symmetric matrices $\bm{\Sigma}_j$ for $j=-1,1$.
We consider the $n$-dimensional two-component Gaussian mixture model for a random vector $\bm{X}$, defined as
\begin{equation}\label{eq:model}
\bm{X} = \theta\bm{\mu} + \bm{g}^{(\theta)}.
\end{equation}

We address the problem of clustering observed random vectors $\bm{X}_1,\ldots,\bm{X}_m$ that are independent and follow the same distribution as $\bm{X}$ in \eqref{eq:model} into two classes determined by the latent sign variable $\theta$.
To this end, we apply spectral clustering algorithm, which classifies each observed random vector according to the sign of $\langle \bm{\gamma}_1(\bm{S}_m), \bm{X}_i \rangle$ for $i=1,\ldots,m$, where $\bm{S}_m   ={ \sum_{i=1}^m \bm {X}_i \bm{X}_i^\top/m}$. 
Since the sign of an eigenvector is indeterminate, we assume that $\langle \bm{\gamma}_1(\bm{S}_m), \bm{\gamma}_1(\bm{\Sigma}) \rangle > 0$, where $\bm{\gamma}_1(\cdot)$ denotes the eigenvector associated with the largest eigenvalue of the matrix argument, and {$\bm{\Sigma} = (\bm{\Sigma}_{-1}+\bm{\Sigma}_1)/2$}.

With this convention, the misclassification probability of $\bm{X}_i$ is represented as
\[{\mathbb P}\big(\,\theta_i \,\langle \bm{\gamma}_1(\bm{S}_m), \bm{X}_i \rangle < 0\,\big)
\] for $i=1,\ldots,m$.


\section{Main result}
Our main result evaluates the misclassification probability of spectral clustering algorithm together with a non-asymptotic error bound.

\begin{thm}\label{mthm3}
Let $c$, $C$, $K$, and $K_g$ be absolute constants independent of $m$, $n$, $\bm{\mu}$, $\bm{\Sigma}_1$, and $\bm{\Sigma}_{-1}$.
Define
\begin{align*}
\delta = & 2^{5/2} c_1 \sqrt{ \frac{n\max_{j=-1,1} \{  \lambda_1(\bm{\Sigma}_j)  \}}{m\| \bm{\mu} \|_2^2}}\\
&+ \sqrt{2}\left(  \frac{\lambda_1(\bm{\Sigma}_1)+\lambda_1(\bm{\Sigma}_{-1})}{\| \bm{\mu} \|_2^2}   \right)
\left( CK^2 \left( \sqrt{\frac{2n}{m}}  + \frac{2n}{m}\right) +1  \right),
\end{align*}
where $c_1 = 1+K_g^2/\sqrt{c}$.
Suppose that there exists $\alpha \in (0,1)$ such that
\begin{align}\label{eq:con}
\delta \leq \frac{\alpha\| \bm{\mu} \|_2}{ c_1\sqrt{n\max_{j=-1,1}   \lambda_1(\bm{\Sigma}_j)} + \| \bm{\mu} \|_2}.
\end{align}
Then, for any $j=-1,1$,
it holds that
\begin{align*}\label{eq:mis_prob_PhiA_refined}
&\left|{\mathbb P}\!\left(j\langle \bm{\gamma}_1(\bm{S}_m),\bm{X}_i\rangle<0\,\middle|\,\theta_i=j\right)-\Phi\skakko{-\frac{\|\bm{\mu}\|_2^2}{\max_{j=-1,1}\|\bm{\Sigma}_j^{1/2}\bm{\mu}\|_2}}\right|\\
&\le \frac{\alpha\|\bm{\mu}\|_2^2}{\max_{j=-1,1}\|\bm{\Sigma}_j^{1/2}\bm{\mu}\|_2}
\,\phi\skakko{\frac{(1-\alpha)\|\bm{\mu}\|_2^2}{\max_{j=-1,1}\|\bm{\Sigma}_j^{1/2}\bm{\mu}\|_2}}+10e^{-n},
\end{align*}
where $\Phi(\cdot)$  and $\phi(\cdot)$ denote the distribution function and density function of $\mathcal{N}(0,1)$, respectively.
\end{thm}

\begin{rem}
Let us briefly explain the absolute constants $c$, $C$, $K$, and $K_g$ in Theorem~\ref{mthm3}. The constants $C$ and $K$ are universal constants that appear in the concentration/perturbation bounds used to control
$\|\bm{S}_m-\bm{\Sigma}\|_{\mathrm{op}}$ and the corresponding eigenvector deviation; in particular, they do not depend on $m$, $n$, $\bm{\mu}$, $\bm{\Sigma}_1$, or $\bm{\Sigma}_{-1}$.
The constant $K_g$ is the sub-Gaussian (Orlicz $\psi_2$) norm of a standard normal random variable. In particular, under the convention $K_g= \inf\left\{ s>0:\ \mathbb{E}[\exp\!\left(Z^2/s^2\right)]\le 2 \right\}$ for
$Z\sim\mathcal{N}(0,1)$, we have $K_g=\sqrt{8/3}$.
Consequently, for $\bm{g}\sim\mathcal{N}_n(\bm{0}_n,\bm{I}_n)$, a standard Gaussian norm concentration inequality yields that for all $t\ge 0$,
\begin{equation}\label{eq:gauss_norm_conc}
\mathbb{P}\bigl(\bigl|\|\bm{g}\|_2-\sqrt{n}\bigr|\ge t\bigr)
\le 2\exp\!\left(-\frac{c\,t^2}{K_g^4}\right),
\end{equation}
see, for example, \citet[Ch.~3]{v18}. The constant $c$ in Theorem~\ref{mthm3} corresponds to the absolute constant in \eqref{eq:gauss_norm_conc}.
With the choice $c_1:=1+K_g^2/\sqrt{c}$, \eqref{eq:gauss_norm_conc} implies
\[
\mathbb{P}\!\left(\|\bm{g}\|_2>c_1\sqrt{n}\right)
\le 2e^{-n},
\]
which is exactly the tail bound used to control the event $\mathrm{E}_2$ in
\eqref{eq:mthm3_eva3}.
\end{rem}

\begin{rem}
When $\bm{\Sigma}_1 = \bm{\Sigma}_2$, 
the quantity 
$\Phi ( -  {\| \bm{\mu} \|_2^2}/{\max_{j=-1,1} \| \bm{\Sigma}_j^{1/2}\bm{\mu} \|_2  } )$ in Theorem \ref{mthm3} reduces to 
$\Phi ( -  {\| \bm{\mu} \|_2^2}/{\| \bm{\Sigma}_j^{1/2}\bm{\mu} \|_2  } )$, whereas the oracle misclassification probability of linear discriminant analysis is given by
$\Phi (-\sqrt{\bm \mu^\top \bm \Sigma_j^{-1}\bm \mu} )$.
We have 
$$
\frac{\| \bm{\mu} \|_2^2}{\left\| \bm{\Sigma}_j^{1/2}\bm{\mu} \right\|_2  }
\leq
\sqrt{\bm \mu^\top \bm \Sigma_j^{-1}\bm \mu}
\leq
\sqrt{\frac{\lambda_1(\bm{\Sigma}_j)}{\lambda_n(\bm{\Sigma}_j)}}
\frac{\| \bm{\mu} \|_2^2}{\left\| \bm{\Sigma}_j^{1/2}\bm{\mu} \right\|_2  },
$$
where $\lambda_n(\bm{\Sigma}_j)$ denotes the smallest eigenvalue of $\bm{\Sigma}_j$,  
which yields
$$
\Phi\skakko{
-\sqrt{\frac{\lambda_1(\bm{\Sigma}_j)}{\lambda_n(\bm{\Sigma}_j)}}
\frac{\| \bm{\mu} \|_2^2}{\left\| \bm{\Sigma}_j^{1/2}\bm{\mu} \right\|_2  }
}
\leq
\Phi\skakko{
-
\sqrt{\bm \mu^\top \bm \Sigma_j^{-1}\bm \mu}
}
\leq
\Phi\skakko{
-
\frac{\| \bm{\mu} \|_2^2}{\left\| \bm{\Sigma}_j^{1/2}\bm{\mu} \right\|_2  }
}.
$$
Thus, our theorem provides a quantitative assessment of the gap between our bound and the oracle misclassification probability.
In the isotropic case $\bm{\Sigma}_1 = \bm{\Sigma}_2= \sigma^2 \bm{I}_n$, the spectral clustering algorithm is optimal.
On the other hand, when $\bm{\Sigma}_1 \neq \bm{\Sigma}_2$, even in an oracle setting, it is difficult to explicitly describe the optimal misclassification probability of linear discriminant analysis \citep{ab62}.
\end{rem}

Let 
\[ \eta :=  \frac{\|\bm\mu\|_2^2}{ \max_{j=-1,1}\lambda_1(\bm\Sigma_j)} \]
that can regarded as the signal-to-noise ratio for our clustering problem.
The numerator of $\eta$ corresponds to the magnitude of separation between the two clusters, while the denominator represents the largest within-cluster variability.
Thus, $\eta$ measures the separation between the cluster means relative to the largest within-cluster dispersion, 
and therefore quantifies the intrinsic difficulty of the clustering problem.
The condition \eqref{eq:con} is satisfied for large $m$ and $n$ if
\[
\frac{1}{\eta}\max\!\left\{
\dfrac{n}{\sqrt{m}},\;
n^{1/3}
\right\}\to 0.
\]

The next corollary demonstrates that spectral clustering algorithm enjoy consistency in a certain high-dimensional scenario.
Intuitively, if $\eta$ is sufficiently large, then the clustering is consistent.

\begin{cor}\label{cor:consistency}
As $n \to \infty$ with $\log m / n \to 0$, if
\[
\frac{1}{\eta}\max\!\left\{
\frac{n}{\sqrt m},\,
n^{1/3},\,
\log m
\right\}\to 0,
\]
then the probability that all observations $\bm X_1,\dots,\bm X_m$ are correctly clustered tends to one.
Moreover, under these conditions, the misclassification probability decays at an exponential rate.
\end{cor}

\begin{rem}
Corollary \ref{cor:consistency} does not require the largest eigenvalue to be bounded.
\end{rem}

Finally, we compare the sufficient conditions for the expected number of misclustered observations to converge to zero, that is
\begin{equation}\label{eq:misc1X}
\mathbb{E} \!\left[
 \min \!\left\{
  \sum_{i=1}^m 1\!\left\{ \theta_i \langle \boldsymbol{\gamma}_1(\boldsymbol{S}_m), \boldsymbol{X}_i \rangle < 0 \right\},
  \sum_{i=1}^m 1\!\left\{ \theta_i \langle \boldsymbol{\gamma}_1(\boldsymbol{S}_m), \boldsymbol{X}_i \rangle > 0 \right\}
 \right\}
 \right] \to 0,
\end{equation}
with those in \cite{cz18} under the isotropic setting $\bm\Sigma_1 = \bm\Sigma_2 = \bm I_n$.
By Corollary~\ref{cor:consistency}, the sufficient conditions for \eqref{eq:misc1X} to hold are
\[
\frac{1}{\|\bm{\mu}\|_2}\max\!\left\{
n^{1/2} m^{-1/4},\;
n^{1/6},\;
(\log m)^{1/2}
\right\}\to0
\text{ and }
\frac{\log m}{n}\to0.
\]
whereas the corresponding conditions in \cite{cz18} are
\[
\frac{1}{\|\bm\mu\|_2}\max\!\left\{
n^{1/4},\; m^{1/2}
\right\} \to0.
\]
When $m$ grows on the same order as the sample size $n$,
our condition is less restrictive. 
In contrast, when $n \to \infty$ while $m$ remains fixed, the conditions in \cite{cz18} is less restrictive.


\section{Proofs}

\subsection{Proof of Theorem~\ref{mthm3}}
\begin{proof}[Proof of Theorem~\ref{mthm3}]
Let
\[
A:=\frac{\|\bm{\mu}\|_2^2}{\max_{j=-1,1}\|\bm{\Sigma}_j^{1/2}\bm{\mu}\|_2}\ (\ge 0)
\quad
\text{and}
\quad
p_j:={\mathbb P}\!\left(j\langle \bm{\gamma}_1(\bm{S}_m),\bm{X}_i\rangle<0\,\middle|\,\theta_i=j\right).
\]
From Lemma~\ref{lemk3} below, it follows that
\[
\Phi\bigl(-(1+\alpha)A\bigr)-10e^{-n}\le p_j\le \Phi\bigl(-(1-\alpha)A\bigr)+10e^{-n}.
\]
Thus, we have
\[ 
|p_j-\Phi(-A)|
\le \max\!\left\{\Phi(-A)-\Phi\bigl(-(1+\alpha)A\bigr),\ \Phi\bigl(-(1-\alpha)A\bigr)-\Phi(-A)\right\}+10e^{-n}.
\]
Since $\Phi'=\phi$ and $\phi(t)=\phi(-t)$ for all $t\in \mathbb{R}$, it holds that 
\[
\Phi(-A)-\Phi\bigl(-(1+\alpha)A\bigr)=\int_A^{(1+\alpha)A}\phi(u)\,du\le \alpha A\,\phi(A),
\]
\[
\Phi\bigl(-(1-\alpha)A\bigr)-\Phi(-A)=\int_{(1-\alpha)A}^{A}\phi(u)\,du\le \alpha A\,\phi\bigl((1-\alpha)A\bigr).
\]
Because $\phi$ is decreasing on $(0,\infty)$ and $(1-\alpha)A\le A$, we have
$\phi(A)\le \phi((1-\alpha)A)$, so the maximum is bounded by $\alpha A\,\phi((1-\alpha)A)$.
\end{proof}

The following lemma provides an bound on the conditional misclassification probability and is key to the proof of Theorem~\ref{mthm3}.
\begin{lem}\label{lemk3}
Let $c$, $C$, and $K$ be absolute constants independent of $m$, $n$, $\bm{\mu}$, $\bm{\Sigma}_1$, and $\bm{\Sigma}_{-1}$.
The quantities $c_1$ and $\delta$ are as defined in Theorem~\ref{mthm3}.
If there exists $\alpha \in (0,1)$ such that
\begin{align}
\delta \leq \frac{\alpha\| \bm{\mu} \|_2}{ c_1\sqrt{n \lambda_1(\bm{\Sigma}_j) } + \| \bm{\mu}\|_2},
\end{align}
it holds that
\begin{align}\label{eq:mthm3_result1}
&\Phi\left( -  \frac{ (1 + \alpha)\| \bm{\mu} \|_2^2}{\left\| \bm{\Sigma}_j^{1/2}\bm{\mu} \right\|_2  } \right) - 10\e^{-n}\\
&\leq
{\mathbb P}(j\langle \bm{\gamma}_1(\bm{S}_m), \bm{X}_i \rangle <0 \arrowvert \theta_i = j )\leq \Phi\left( -  \frac{ (1 - \alpha)\| \bm{\mu} \|_2^2}{\left\| \bm{\Sigma}_j^{1/2}\bm{\mu} \right\|_2  } \right) + 10\e^{-n}
\end{align}
for $i=1,\ldots,m$ and $j=-1,1$.
\end{lem}

\begin{proof}
We prove \eqref{eq:mthm3_result1} only for $j=1$, and the case $j=-1$ follows similarly.
For notational convenience, define the events $\mathrm{A}_i$, $\mathrm{E}_1$, and $\mathrm{E}_2$ as
\begin{align*}
\mathrm{A}_i &= \left\{ \theta_i = 1 \right\} ,\quad
\mathrm{E}_1 = \left\{ \|\bm{\gamma}_1(\bm{S}_m) - \tilde{\bm{\mu}} \|_2 \leq \delta   \right\},  \\
\mathrm{E}_2 &= \left\{  \left\| \bm{g}_i \right\|_2  \leq c_1\sqrt{n\lambda_1(\bm{\Sigma}_1)} \right\},
\end{align*}
where $\tilde{\bm{\mu}} = \bm{\mu}/\|\bm{\mu}\|_2$.
Then, we find that
\begin{align}
&{\mathbb P}\left(  \left\langle \bm{\gamma}_1(\bm{S}_m), \bm{X}_i \right\rangle < 0 \middle\arrowvert \mathrm{A}_i \right)\nonumber \\
\leq&
{\mathbb P}\left( \left\{ \left\langle \bm{\gamma}_1(\bm{S}_m), \bm{X}_i \right\rangle < 0 \right\} \cap \mathrm{E}_1 \cap \mathrm{E}_2 \middle\arrowvert \mathrm{A}_i \right)
 + {\mathbb P}\left(  \mathrm{E}_1^c \middle\arrowvert \mathrm{A}_i \right) +  {\mathbb P}\left(  \mathrm{E}_2^c \middle\arrowvert \mathrm{A}_i \right). \label{eq:mthm3_eva1}
\end{align}

First, we evaluate the first term on the right-hand side of \eqref{eq:mthm3_eva1}:
\begin{align}
&{\mathbb P} \biggl(\langle \bm{\gamma}_1(\bm{S}_{m}) , \bm{X}_i \rangle<0 \biggl\arrowvert \mathrm{A}_i\biggr) \nonumber \\*
=& {\mathbb P} \biggl(\langle \bm{\gamma}_1(\bm{S}_{m}) - \tilde{\bm{\mu}} + \tilde{\bm{\mu}} , \theta_i \bm{\mu} + \bm{g}_i^{(\theta_i)}  \rangle<0 \biggl\arrowvert \mathrm{A}_i \biggr) \nonumber \\*
=& {\mathbb P}\biggl(  \langle \bm{\gamma}_1(\bm{S}_{m}) - \tilde{\bm{\mu}},  \bm{\mu} \rangle + \langle \bm{\gamma}_1(\bm{S}_{m}) -\tilde{\bm{\mu}},\bm{g}_i^{(\theta_i)}\rangle +\langle          \tilde{\bm{\mu}} ,\bm{\mu} \rangle \nonumber \\&\quad +\langle\tilde{\bm{\mu}},\bm{g}_i^{(\theta_i)}\rangle <0 \biggl\arrowvert \mathrm{A}_i\biggr) \nonumber \\*
\leq & 
{\mathbb P}\biggl( - \|\bm{\gamma}_1(\bm{S}_{m}) - \tilde{\bm{\mu}}\|_2\cdot \| \bm{\mu} \|_2  + \langle\bm{\Sigma}_1^{1/2}(\bm{\gamma}_1(\bm{S}_{m})- \tilde{\bm{\mu}}),\bm{\Sigma}_1^{-1/2}\bm{g}_i^{(\theta_i)}\rangle\nonumber \\*
&\qquad + \| \bm{\mu}\|_2  +  \|\bm{\Sigma}_1^{1/2}\tilde{\bm{\mu}}\|_2 \biggl\langle \frac{\bm{\Sigma}_1^{1/2}\tilde{\bm{\mu}}}{\| \bm{\Sigma}_1^{1/2}\tilde{\bm{\mu}}\|_2 },\bm{\Sigma}_1^{-1/2}\bm{g}_i^{(\theta_i)}\biggr\rangle
 < 0\biggl\arrowvert \mathrm{A}_i \biggr) \nonumber \\*
\leq &
{\mathbb P}\biggl(
 - \|\bm{\gamma}_1(\bm{S}_{m}) - \tilde{\bm{\mu}}\|_2\cdot \| \bm{\mu} \|_2 -\|\bm{\Sigma}_1^{1/2} \|_{\mathrm{op}} \|\bm{\gamma}_1(\bm{S}_{m}) - \bm{\gamma}_1 (\bm{\Sigma}_1)\|_2 \nonumber\\ 
&\cdot \|\bm{\Sigma}_1^{- 1/2}\bm{g}_i^{(\theta_i)}\|_2
+ \| \bm{\mu}\|_2 +  \|\bm{\Sigma}_1^{1/2}\tilde{\bm{\mu}}\|_2 \biggl\langle \frac{\bm{\Sigma}_1^{1/2}\tilde{\bm{\mu}}}{\| \bm{\Sigma}_1^{1/2}\tilde{\bm{\mu}}\|_2 },\bm{\Sigma}_1^{-1/2}\bm{g}_i^{(\theta_i)}\biggr\rangle
 < 0\biggl\arrowvert \mathrm{A}_i \biggr) \label{eq:mthm3_f1}.
\end{align} 
Here, we used the following inequality
\begin{align*}
&\langle \bm{\gamma}_1(\bm{S}_{m}) - \bm{\gamma}_1 (\bm{\Sigma}_1), \bm{g}_i^{(\theta_i)} \rangle  
= \langle \bm{\Sigma}_1^{1/2} (\bm{\gamma}_1(\bm{S}_{m}) - \bm{\gamma}_1 (\bm{\Sigma}_1)), \bm{\Sigma}_1^{-1/2}\bm{g}_i^{(\theta_i)}\rangle  \\
&\geq -
   \|\bm{\Sigma}_1^{1/2} (\bm{\gamma}_1(\bm{S}_{m}) - \bm{\gamma}_1 (\bm{\Sigma}_1))\|_2 \|\bm{\Sigma}_1^{- 1/2}\bm{g}_i^{(\theta_i)}\|_2   \\
&\geq -\|\bm{\Sigma}_1^{1/2} \|_{\mathrm{op}} \|\bm{\gamma}_1(\bm{S}_{m}) - \bm{\gamma}_1 (\bm{\Sigma}_1)\|_2  \|\bm{\Sigma}_1^{- 1/2}\bm{g}_i^{(\theta_i)}\|_2 .
\end{align*}
Set
\begin{align*}
\bm{g}_i  &= \bm{\Sigma}_1^{-1/2}\bm{g}_i^{(\theta_i)}  \quad\text{and}\quad
Z_i =  \biggl\langle \frac{\bm{\Sigma}_1^{1/2}\tilde{\bm{\mu}}}{\|\bm{\Sigma}_1^{1/2}\tilde{\bm{\mu}}\|_2 },\bm{g}_i \biggr\rangle .
\end{align*}
Because $\bm{g}_i \mid \mathrm{A}_i \sim \mathcal{N}_n(0,I_n)$, for any deterministic $a\in\mathbb{R}^n$ we have
$\bm{a}^\top \bm{g}_i \mid \mathrm{A}_i \sim \mathcal{N}(0,\|\bm{a}\|_2^2)$.
By taking
\[
a=\frac{\bm{\Sigma}_1^{1/2}\tilde{\bm{\mu}}}{\bigl\|\bm{\Sigma}_1^{1/2}\tilde{\bm{\mu}}\bigr\|_2},
\]
it follows that
\begin{equation*}
\|\bm{\Sigma}_1^{1/2}\tfrac{\bm{\mu}}{\| \bm{\mu} \|_2}\|_2\, Z_i \mid \mathrm{A}_i \sim \mathcal{N}\!\left(0, \bigl\|\bm{\Sigma}_1^{1/2}\tfrac{\bm{\mu}}{\| \bm{\mu} \|_2}\bigr\|_2^2\right).
\end{equation*}
Hence, the right-hand side of \eqref{eq:mthm3_f1} equals
\begin{align*}
&{\mathbb P}\Biggl( - \| \bm{\gamma}_1(\bm{S}_{m}) - \tilde{\bm{\mu}}\|_2 \left(\sqrt{\lambda_1(\bm{\Sigma}_1)}\| \bm{g}_i\|_2 + \|\bm{\mu} \|_2\right)  + \|\bm{\mu} \|_2 \\
&\quad +  \|\bm{\Sigma}_1^{1/2}\tilde{\bm{\mu}}\|_2 Z_i <0\biggl\arrowvert \mathrm{A}_i \Biggr)\\
=&
{\mathbb P}\Biggl(\|\bm{\Sigma}_1^{1/2}\tilde{\bm{\mu}}\|_2Z_i<  \| \bm{\gamma}_1(\bm{S}_{m}) - \tilde{\bm{\mu}}\|_2  \left(\sqrt{\lambda_1(\bm{\Sigma}_1)}\| \bm{g}_i \|_2 + \| \bm{\mu} \|_2 \right) \\
&\quad - \| \bm{\mu} \|_2\biggl\arrowvert \mathrm{A}_i \Biggr).
\end{align*}
Since
$
\delta \leq {\alpha\| \bm{\mu} \|_2}/{(c_1\sqrt{n\  \lambda_1(\bm{\Sigma}_1) } + \| \bm{\mu} \|_2)}$,
it holds that
\begin{align}
&{\mathbb P}\Biggl( \Biggl
\{\|\bm{\Sigma}_1^{1/2}\tilde{\bm{\mu}}\|_2Z_i<  \| \bm{\gamma}_1(\bm{S}_{m}) - \tilde{\bm{\mu}}\|_2 \left( \sqrt{\lambda_1(\bm{\Sigma}_1)}\| \bm{g}_i \|_2 + \| \bm{\mu} \|_2 \right) - \| \bm{\mu} \|_2 \Biggr\} \nonumber\\
& \quad \cap  \mathrm{E}_1 \cap   \mathrm{E}_2\biggl\arrowvert \mathrm{A}_i \Biggr) \nonumber\\
\leq&  {\mathbb P}\Biggl( \Biggl\{\|\bm{\Sigma}_1^{1/2}\tilde{\bm{\mu}}\|_2Z_i<  \delta \left(c_1\sqrt{n\lambda_1(\bm{\Sigma}_1)} + \| \bm{\mu} \|_2 \right) - \| \bm{\mu} \|_2 \Biggr\} \cap  \mathrm{E}_1 \cap   \mathrm{E}_2 \biggl\arrowvert \mathrm{A}_i \Biggr) \nonumber\\
\leq&
{\mathbb P} \Biggl(Z_i< \frac{-(1-\alpha) \| \bm{\mu} \|_2}{\|\bm{\Sigma}_1^{1/2}\tilde{\bm{\mu}}\|_2} \biggl\arrowvert \mathrm{A}_i \Biggr) \nonumber\\
=&
{\mathbb P} \Biggl(Z_i< \frac{-(1-\alpha) \| \bm{\mu} \|_2^2}{\|\bm{\Sigma}_1^{1/2}\bm{\mu}\|_2} \biggl\arrowvert \mathrm{A}_i \Biggr) .\label{eq:mthm3_eva2}
\end{align}

Next, we evaluate the remaining two terms in \eqref{eq:mthm3_eva1}.
Since $c_1  = 1 + { K_g^2}/{\sqrt{c}}$ by assumption, a standard concentration bound for Gaussian norms gives
\begin{align}
{\mathbb P}(\|\bm{g}_i\|_2 > c_1 \sqrt{n}  \arrowvert \mathrm{A}_i ) 
=& {\mathbb P}( |  \| \bm{g}_i \|_2 - \sqrt{n}  + \sqrt{n} | > c_1 \sqrt{n} \arrowvert \mathrm{A}_i ) \nonumber \\
 \leq& {\mathbb P}( |  \| \bm{g}_i \|_2 - \sqrt{n}  | > (c_1-1) \sqrt{n} \arrowvert \mathrm{A}_i ) \nonumber \\
 \leq& 2 \e^{-n}, \label{eq:mthm3_eva3}
\end{align}
see, e.g., \citet[Ch.~3]{v18}. 
We also find that
\begin{align}
{\mathbb P}\left(\left\|\bm{\gamma}_1(\bm{S}_{m}) - \tilde{\bm{\mu}} \right\|_2 > \delta  \biggl\arrowvert \mathrm{A}_i  \right)\leq&
\frac{1}{{\mathbb P}\left( \mathrm{A}_i \right)}{\mathbb P}\left(\left\|\bm{\gamma}_1(\bm{S}_{m}) - \tilde{\bm{\mu}} \right\|_2 > \delta    \right) \nonumber \\
\leq&
8\e^{-n}. \label{eq:mthm3_eva4}
\end{align}
Combining \eqref{eq:mthm3_eva1}, \eqref{eq:mthm3_eva2}, \eqref{eq:mthm3_eva3}, and \eqref{eq:mthm3_eva4} yields
\begin{equation*}
{\mathbb P} \biggl(\langle \bm{\gamma}_1(\bm{S}_{m}), \bm{X}_i\rangle<0 \biggl\arrowvert \mathrm{A}_i \biggr) \\
\leq {\mathbb P} \Biggl(Z_i< \frac{-(1-\alpha) \| \bm{\mu} \|_2^2}{\|\bm{\Sigma}_1^{1/2}\bm{\mu}\|_2} \biggl\arrowvert \mathrm{A}_i \Biggr) +  10 \e^{-n}.
\end{equation*}
Similarly, by an analogous argument, we obtain
\begin{align*}
&{\mathbb P}\bigl(\langle \bm{\gamma}_1(\bm{S}_{m}) , \bm{X}_i \rangle<0 \mid \mathrm{A}_i\bigr) \\
&\geq {\mathbb P}\Bigl(
  \bigl\|\bm{\gamma}_1(\bm{S}_{m}) - \tfrac{\bm{\mu}}{\| \bm{\mu} \|_2}\bigr\|_2
  \bigl(\sqrt{\lambda_1(\bm{\Sigma}_1)}\| \bm{g}_i\|_2 + \|\bm{\mu} \|_2\bigr)
  + \|\bm{\mu} \|_2
  +  \bigl\|\bm{\Sigma}_1^{1/2}\tfrac{\bm{\mu}}{\| \bm{\mu} \|_2}\bigr\|_2 Z_i <0
  \,\bigm|\, \mathrm{A}_i
\Bigr)\\
&\geq {\mathbb P}\Bigl(
  \bigl\{
    \bigl\|\bm{\Sigma}_1^{1/2}\tfrac{\bm{\mu}}{\| \bm{\mu} \|_2}\bigr\|_2 Z_i
      < - \bigl\|\bm{\gamma}_1(\bm{S}_{m}) - \tfrac{\bm{\mu}}{\| \bm{\mu} \|_2}\bigr\|_2
        \bigl( \sqrt{\lambda_1(\bm{\Sigma}_1)}\| \bm{g}_i \|_2 + \| \bm{\mu} \|_2 \bigr)
        - \| \bm{\mu} \|_2
  \bigr\} \\
&\qquad\qquad
  \cap  \mathrm{E}_1 \cap   \mathrm{E}_2
  \,\bigm|\, \mathrm{A}_i
\Bigr)\\
&\geq {\mathbb P}\Bigl(
  \bigl\{
    \bigl\|\bm{\Sigma}_1^{1/2}\tilde{\bm{\mu}}\bigr\|_2 Z_i
      < - \delta \bigl(c_1\sqrt{n\lambda_1(\bm{\Sigma}_1)} + \| \bm{\mu} \|_2 \bigr)
        - \| \bm{\mu} \|_2
  \bigr\}
  \cap  \mathrm{E}_1 \cap   \mathrm{E}_2
  \,\bigm|\, \mathrm{A}_i
\Bigr)\\
&\geq
{\mathbb P}\Biggl(
  Z_i <
  \frac{-(1+\alpha)\,\|\bm{\mu}\|_2}
       {\bigl\|\bm{\Sigma}_1^{1/2}\tilde{\bm{\mu}}\bigr\|_2}
  \,\biggm|\, \mathrm{A}_i
\Biggr)
- {\mathbb P}(\mathrm{E}_1^c\mid \mathrm{A}_i)
- {\mathbb P}(\mathrm{E}_2^c\mid \mathrm{A}_i)\\
&\geq
{\mathbb P}\Biggl(
  Z_i <
  \frac{-(1+\alpha)\,\|\bm{\mu}\|_2}
       {\bigl\|\bm{\Sigma}_1^{1/2}\tilde{\bm{\mu}}\bigr\|_2}
  \,\biggm|\, \mathrm{A}_i
\Biggr)
-10\e^{-n}.
\end{align*}
\end{proof}



\begin{funding}
This study was supported in part by JSPS Grant-in-Aid for Early-Career Scientists JP23K16851 (Y.G.), JSPS KAKENHI Grant Number JP25K07133 (K.T.), and the Research Fellowship Promoting International Collaboration of the Mathematical Society of Japan (Y.G.).
\end{funding}

\bibliographystyle{imsart-nameyear}
\bibliography{ref}

\end{document}